\documentclass[12pt ]{article}

\newcommand{\bc}{\begin{center}}
\newcommand{\ec}{\end{center}}
\newcommand{\ba}{\begin{array}}
\newcommand{\ea}{\end{array}}

\date{}

\newcommand{\na}{\nabla}

\newcommand{\wna}{\tilde{\nabla}}
\newcommand{\yrd}{1 + \frac{d(x,\,x_{0})}{r}}
\newcommand{\gij}{g_{ij}(x)}
\newcommand{\wgij}{\tilde{g}_{ij}(x)}

\newcommand{\vob}{Vol(B(x_{0},\,r))}
\newcommand{\vr}{V_{0}(r)}
\newcommand{\cm}{C_{0}^{\infty}(M)}
\newcommand{\rn}{R^{n}}

\newcommand{\py}{\frac{1}{p}}
\newcommand{\qy}{\frac{1}{q}}
\newcommand{\ny}{\frac{1}{n}}
\newcommand{\nc}{\tilde{C}}

\title{ Manifolds with
 non-negative Ricci curvature and Nash inequalities
 \footnotetext{ Supported by NSFC  No:
 10271089;}}
\author{Qihua Ruan and Zhihua Chen}

\begin{document}

\maketitle
\thispagestyle{empty}


{\bf Abstract.} We prove that for any complete n-dimensional
Riemannian manifold with nonnegative Ricci curvature, if the Nash
inequality is satisfied, then it is diffeomrphic to $R^{n}$



\setcounter{equation}{0} 
$$\ $$

{\bf 1. Introduction. }Let M be any complete n-dimensional ($n
\geq 2$ ) Riemannian manifold with nonnegative Ricci curvature,
$\cm$ be the space of smooth functions with compact support in M.
Denote by dv and $\na$ the Riemannian volume element and the
gradient operator of M, respectively.

It is well known in [1] that Ledoux showed that:  If one of the
following Sobolev inequalities is satisfied,
$$||f||_{p}\leq
C_{0}||\na f||_{q},\forall f \in \cm,\,\, 1 \leq q <n,\,\,\py =
\qy - \ny.\leqno(1)$$ where $C_{0} $ is the optimal constant in
$\rn,$  denote $||f||_{p}$ by the $ L^{p}$ norm of function $f$
;then M is \ isometric to $\rn.$

The basic idea of Ledoux's result is to find a function in $\cm$,
then one can substitute it to (1) and obtain that $\vob
\geq\vr$,here $\vob$ denote the volume of the geodesic ball
$B(x_{0},r)$ of radius $r$ with center $x_{0}$,and $\vr$ the
volume of the Euclidean ball of radius $r$ in $\rn.$ Since the
Ricci curvature of M is nonnagtive ,from Bishop's comparison
theorem[2],we know that $\vob \leq \vr,$ so M is isometric to
$\rn.$

Later Xia combined Ledoux's method with Cheeger and Colding's
result[3],which is that given an integer $n\geq 2$,there exists a
constant $\delta (n)>0$ such that any n-dimensional complete
Riemannian manifold with nonnegative Ricci curvature and $\vob
\geq(1-\delta(n))\vr$ for some $x_{0}\in M$ and all $ r>0 $ is
diffeomorphic to $\rn$. He proved that [4]: If one of the
following sobolev inequalities is satisfied,
$$||f||_{p}\leq
C_{1}||\na f||_{q},\forall f \in \cm,\,\, 1 \leq q <n,\,\,\py =
\qy - \ny.\leqno(2)$$ where the positive constant
$C_{1}>C_{0}$,then M is diffiomorphic to $\rn$.

This is a beautiful result.As we known,the Sobolev inequality(2)
belongs to a general family of inequalities of the type
$$||f||_{p}\leq
C||f||_{s}^{\theta}\ ||\na f||_{q}^{1-\theta},\forall f \in
\cm,\,\,
\frac{1}{t}=\frac{\theta}{s}+\frac{1-\theta}{p}.\leqno(3)$$ (see
[5]).Inequality (2) corresponds to $\theta = 0$. When $q=r=2$,and
$\theta=2/(n+2)$,it corresponds to the Nash inequality
$$\left(\int\mid f \mid^{2}dv\right)^{1+\frac{2}{n}} \leq C\left(\int\mid f \mid dv\right)
^{\frac{4}{n}} \int\mid \nabla f \mid^{2}dv,\,f \in
C_{0}^{\infty}(M),\leqno(4)$$ [ see(6)].So we may naturally ask
whether or not there has analogous result for the Nash inequality?
In this note ,we confirmed this problem.

{\sl {MAIN THEOREM } }{\sl Let M be a complete n-dimensional
Riemannian manifold with non-negative Ricci curvature. If the Nash
inequality $(4)$ is satisfied with the positive constant C, then M
is diffeomorphic to $R^{n}.$ }

\setcounter{equation}{0}
{\bf 2. Proof of Main Theorem.} Before showing this theorem, we
must introduce a lemma about Scheon and Yau's cut-off function
(see[7]or[8]), because we will use it.

{\sl LEMMA 2.1 } Suppose $(M,\,\gij)$ is an n-dimensional complete
Riemannian manifold with non-negative Ricci curvature. Then there
exists a constant $\nc$ depending only on the dimension n, such
that for any $x_{0} \in M$ and any number $0 < r < +\infty,$ there
exists a smooth function $\varphi(x) \in C^{\infty}(M)$ satisfying
$$\left\{\ba{ll}e^{\nc (\yrd)} &\leq \varphi(x) \leq e^{-(\yrd)}\\\,\,\,\,\,\,\,\,
 \mid \na \varphi(x) \mid &\leq \frac{\nc}{r}\varphi(x) \ea \right.\leqno(5)$$
 for $\forall x \in M,$ where $d(x,\,x_{0})$ denote the distance between x and
 $x_{0}$ with respect to the metric $\gij.$

{\sl PROOF.}\ \ In [7], Scheon and Yau constructed a $C^{\infty}$
function $\psi(x)$ satisfying $$\left\{\ba{ll}\frac{1}{\nc}(1 +
d(x,\,x_{0})) &\leq \psi(x) \leq 1 +
d(x,\,x_{0})\\\,\,\,\,\,\,\,\,
 \mid \na \psi(x) \mid &\leq \nc \ea \right.$$
 for $\forall x \in M$ and some positive constant $\nc$ depending only on the
 dimension n. Now we define a new metric on M by $$\wgij =
 \frac{1}{r^{2}}\gij,\,\,\,\,x \in M.$$ Then the new metric $\wgij$ is still a
 complete Riemannian metric on M with non-negative Ricci curvature. Thus there
 exists a smooth function $h(x) \in C^{\infty}(M)$ such that
 $$\left\{\ba{ll}\frac{1}{\nc}(\yrd) &\leq h(x) \leq -(\yrd)\\\,\,\,\,\,\,\,\,
 \mid \wna h(x) \mid_{\wgij} &\leq \nc \ea\right.$$
for $\forall x \in M.$ Here $\wna$ and $\mid \cdot \mid_{\wgij}$
are the gradient operator and norm with respect to the new metric
$\wgij.$ Then by setting $\varphi(x) = e^{-h(x)}$ we get the
desired cut-off function (4).

As following Ledoux's method,we want to look for a function in
$\cm.$ we just obtain a smooth function from Lemma$ \  2.1$ ,if we
substitute it to (4),then through direct compute we can get that
$\vob\geq C \vr \  \ (0<C<1)$. From Cheeger and Colding's result
,we can prove the Main Theorem. Now the question is that
$\varphi(x)$ hasn't compact support in M. To solve this problem,
we can find a sequence functions $\varphi_{m}(x) \in \cm$ such
that $||\varphi_{m}||_{t}\rightarrow ||\varphi||_{t}$ and
$||\nabla\varphi_{m}||_{t} \rightarrow ||\nabla\varphi||_{t}$,
when $m \rightarrow +\infty,$for any $t>0.$ Then from Lebesgue
dominated convergence theorem,we get $\varphi(x)$ is satisfied
with the inequality (4).


{\sl PROPOSITION 2.1}  $ \ \ $For $\varphi(x)$ in (4), there exist
a sequence functions $\varphi_{m}(x) \in \cm$ such that
$||\varphi_{m}||_{t}\rightarrow ||\varphi||_{t}$ and
$||\nabla\varphi_{m}||_{t} \rightarrow ||\nabla\varphi||_{t}$,
when $m \rightarrow +\infty,$for any $t>0.$

{\sl PROOF.}\ \ As we known, there always exists functions
$\psi_{m}(x) \in \cm$ such that $\psi_{m}(x)=1 $ for $ x \in
B(x_{0},2^{m}r)$,any fix positive number r, $\psi_{m}(x)=0 $ for $
x \in M\setminus B(x_{0},2^{m+1}r)$, otherwise $0\leq
\psi_{m}(x)\leq 1$;and $|\nabla \psi_{m}|\leq \frac{2}{2^{m}r}.$
Let $\varphi_{m}(x)=\psi_{m}(x)\varphi(x),$ then $\varphi_{m}(x)
\in \cm$.Thus
$$||\varphi_{m}-\varphi||_{t}^{t}=\int_{M}|\psi_{m}\varphi-\varphi|^{t}dv
=\int_{M\setminus B(x_{0},2^{m}r)}|\psi_{m}\varphi-\varphi|^{t}dv
\leq \int_{M\setminus B(x_{0},2^{m}r)}|\varphi|^{t}dv $$ then we
only need to prove that$$ \int_{M\setminus
B(x_{0},2^{m}r)}|\varphi|^{t}dv\longrightarrow 0, when\
m\longrightarrow
 +\infty.\leqno (5)$$
 From Lemma $2.1$,we have $$\int_{M\setminus B(x_{0},2^{m}r)}|\varphi|^{t}dv
 \leq \int_{M\setminus B(x_{0},2^{m}r)}exp (
 -t(1+\frac{d(x_{0},x)}{r})dv$$
$$\leq C
\int_{2^{m}r}^{+\infty}e^{-\frac{t}{r}\rho}\rho^{n-1}d\rho$$
$$\leq -C\rho^{n-1}e^{-\frac{t}{r}\rho}\mid_{2^{m}r}^{+\infty}+C
\int_{2^{m}r}^{+\infty}e^{-\frac{t}{r}\rho}\rho^{n-2}d\rho,
(m\longrightarrow
 +\infty)$$
 $$\leq C
 \int_{2^{m}r}^{+\infty}e^{-\frac{t}{r}\rho}\rho^{n-2}d\rho$$
 $$\leq...\leq C\int_{2^{m}r}^{+\infty}e^{-\frac{t}{r}\rho}d\rho$$
 $$\leq -Ce^{-\frac{t}{r}\rho}\mid_{2^{m}r}^{+\infty}\longrightarrow
 0,(m \longrightarrow +\infty).$$
So we obtain that $||\varphi_{m}||_{t}\rightarrow ||\varphi||_{t},
m \longrightarrow +\infty.$ From (4) and $|\nabla\psi_{m}|\leq
\frac{2}{2^{m}r}$, we get that
$$||\nabla\varphi_{m}-\nabla\varphi||_{t}=||\nabla\psi_{m}\varphi+\psi_{m}\nabla\varphi-\nabla\varphi||_{t}$$
$$\leq||\nabla\psi_{m}\varphi||_{t}+||(\psi_{m}-1)\nabla\varphi||_{t}$$
$$\leq(\int_{M\setminus B(x_{0},2^{m}r)}|\nabla\psi_{m}|^{t}|\varphi|^{t}dv)^{\frac{1}{t}}+
(\int_{M\setminus
B(x_{0},2^{m}r)}|\psi_{m}-1|^{t}|\nabla\varphi|^{t}dv)^{\frac{1}{t}}$$
$$\leq\frac{2}{2^{m}r}(\int_{M\setminus
B(x_{0},2^{m}r)}|\varphi|^{t}dv)^{\frac{1}{t}} +
\frac{\nc}{r}(\int_{M\setminus
B(x_{0},2^{m}r)}|\varphi|^{t}dv)^{\frac{1}{t}}$$ Which combining
with (5) implies that  $||\nabla\varphi_{m}||_{t} \rightarrow
||\nabla\varphi||_{t}$,
 $m \rightarrow +\infty.$

 Now we prove the Main Theorem.

 {\sl PROOF:}\ \ By Proposition $2.1$ we know the function
 $\varphi(x)$ satisfies with the Nash inequality.Together with (4)
 we get
$$
\left(\int \mid \varphi \mid^{2}dv \right)^{1+\frac{2}{n}} \leq
C\left(\int \mid \varphi \mid dv \right)^{\frac{4}{n}}\left(\int
\mid \nabla \varphi \mid^{2}dv
\right)\,\,\,\,\,\,\,\,\,\,\,\,\,\,\,$$
 $$ \left(\int \mid
\varphi \mid^{2}dv \right)^{\frac{2}{n}} \leq C\left(\int \mid
\varphi \mid dv
\right)^{\frac{4}{n}}\left(\frac{\nc^{2}}{r^{2}}\int \mid \varphi
\mid^{2}dv \right)$$
 $$ \left(\int \mid \varphi \mid^{2}dv
\right)^{\frac{2}{n}} \leq \frac{C \nc^{2}}{r^{2}}\left(\int \mid
\varphi \mid dv
\right)^{\frac{4}{n}}\,\,\,\,\,\,\,\,\,\,\,\,\,\,\,\,\,\,\,\,\,\,\,\,\,\,\,\,\,\,\,\,\,\,\,\,\,\,\,\,\,$$
$$\left(\int _{M}\mid \varphi \mid^{2}dv \right) \leq
\left(\frac{C \nc^{2}}{r^{2}}\right)^{\frac{n}{2}}\left(\int_{M}
\mid \varphi \mid dv \right)^{2}
\,\,\,\,\,\,\,\,\,\,\,\,\,\,\,\,\,\,\,\,$$
$$\int_{M}e^{-2\nc\left(1+\frac{d(x,\,x_{0})}{r}\right)}dv \leq
\left(\frac{C \nc^{2}}{r^{2}}\right)^{\frac{n}{2}}
\left(\int_{M}e^{-\left(1+\frac{d(x,\,x_{0})}{r}\right)}dv\right)^{2}\,\,\,\,\,\,\,\,\,\,$$
$$\int_{B(x_{0},\,r)}e^{-2\nc\left(1+\frac{d(x,\,x_{0})}{r}\right)}
\leq \left(\frac{C
\nc^{2}}{r^{2}}\right)^{\frac{n}{2}}\left[vol(B(x_{0},\,r))
+\sum_{k=0}^{+\infty}\int_{B(x_{0},\,2^{k+1}r)\setminus
B(x_{0},\,2^{k}r)}
e^{-2\left(1+\frac{d(x,\,x_{0})}{r}\right)}\right]^{2}$$
$$e^{-4\nc}vol(B(x_{0},\,r)) \leq \left(\frac{C
\nc^{2}}{r^{2}}\right)^{\frac{n}{2}}\left[vol(B(x_{0},\,r)) +
\sum_{k=0}^{+\infty}e^{-2^{k}}(2^{k+1})^{2n}vol(B(x_{0},\,r))\right]^{2}$$
$$e^{-4\nc}vol(B(x_{0},\,r)) \leq \left(\frac{C
\nc^{2}}{r^{2}}\right)^{\frac{n}{2}}C^{2}_{2}(vol(B(x_{0},\,r)))^{2}$$
$$vol(B(x_{0},\,r)) \geq e^{-4\nc}(C
\nc^{2})^{-\frac{n}{2}}C_{2}^{-2}r^{n}$$

 Let $C =e^{-4\nc}(C
\nc^{2})^{-\frac{n}{2}}C_{2}^{-2},$ where $C_{2} = 1 +
\sum_{k=0}^{+\infty}e^{-2^{k}}(2^{k+1})^{2n},$ then
$vol(B(x_{0},\,r)) \geq C r^{n}.$

  Then from above discussion,we know that there exists a number $\delta (n) $,($0<\delta (n)<1 $),such that
  $vol(B(x_{0},\,r))\geq (1-\delta (n))\vr $,together with Cheeger and Colding's result, which implies
  M is diffeomorphic
 to $\rn.$\ \ The proof of the Main Theorem is completed.

\end{document}